\documentstyle[fleqn,12pt]{article}
\begin{document}
\pagenumbering{arabic}
\setcounter{page}{1}
\pagestyle{plain}
\baselineskip=16pt

\thispagestyle{empty}
\rightline{YTUMB 2001-02, November 2001} 
\vspace{1.4cm}

\begin{center}
{\Large\bf Differential geometry of the Z$_3$-graded \\ quantum superplane}
\end{center}

\vspace{1cm}
\begin{center} Salih \c Celik \footnote{E-mail: sacelik@yildiz.edu.tr}
 
Yildiz Technical University, Department of Mathematics, \\
34210 DAVUTPASA-Esenler, Istanbul, TURKEY. \end{center}

\vspace{2cm}
{\bf Abstract}

In this work, differential geometry of the Z$_3$-graded quantum superplane 
is constructed. The corresponding quantum Lie superalgebra and its Hopf 
algebra structure are obtained. 

\vfill\eject
\noindent
{\bf 1. Introduction}

\noindent
Noncommutative geometry [1] has started to play an important role 
in different fields of mathematical physics over the past decade. 
The basic structure giving a direction to the noncommutative geometry 
is a differential calculus on an associative algebra. The noncommutative 
differential geometry of quantum groups was introduced by Woronowicz [2]. 
In this approach the quantum group is taken as the basic noncommutative 
space and the differential calculus on the group is deduced from the 
properties of the group. The other approach, initiated by Wess and Zumino [3], 
followed Manin's emphasis [4] on the quantum spaces as the primary objects. 
Differential forms are defined in terms of noncommuting (quantum) coordinates, 
and the differential and algebraic properties of quantum groups acting on 
these spaces are obtained from the properties of the spaces. The natural 
extension of their scheme to superspace [5] was introduced in [6,7]. 

Recently, there have been many attempts to generalize Z$_2$-graded 
constructions to the Z$_3$-graded case [8-12]. Chung [12] studied the 
Z$_3$-graded quantum space that generalizes the Z$_2$-graded space called a 
superspace, using the methods of Wess and Zumino [3]. In this work, we have 
investigated the noncommutative geometry of the Z$_3$-graded quantum 
superplane. These calculi are discussed from the covariance point of view, 
using the Hopf algebra structure of the quantum superplane [13]. In order to 
obtain the corresponding quantum Z$_3$-grading Lie superalgebra, we 
constructed a left-covariant differential calculus on the Z$_3$-graded quantum 
superplane (of course, this may be done using a right-covariant differential 
calculus on it, too). Hopf algebra structure of the obtained superalgebra is 
given, using the method in [14]. 

Let us briefly investigate a general Z$_3$-graded algebraic structure. Let $z$ 
be a Z$_3$-graded  variable. Then we say that the variable $z$ satisfies the 
relation 
$$z^3 = 0.$$
If $f(z)$ is an arbitrary function of the variable $z$, then the function 
$f(z)$ becomes a polynomial of degree two in $z$, that is, 
$$f(z) = a_0 + a_1 z + a_2 z^2,$$
where $a_0$, $a_2$, $a_1$ denote three fixed numbers whose grades are 
$grad (a_0) = 0$, $grad (a_2) = 1$ and $grad(a_1) = 2$, respectively. 

The cyclic  group Z$_3$ can be represented in the complex plane by means of 
the cubic roots of 1: let $j = e^{{2\pi i}\over 3}$ $(i^2 = - 1)$. Then one 
has 
$$j^3 = 1 \quad \mbox{and} \quad j^2 + j +1 = 0, \quad \mbox{or} \quad 
  (j + 1)^2 = j. $$
One can define the Z$_3$-graded commutator $[A,B]$ as 
$$[A,B]_{Z_3} = AB - j^{ab} BA, $$
where $grad(A) = a$ and $grad(B) = b$. If $A$ and $B$ are $j$-commutative, 
then we have 
$$AB = j^{ab} BA.$$

\vspace*{0.3cm}\noindent
{\bf 2. The algebra of functions on the Z$_3$-graded quantum superplane} 

\noindent
It is well known that the Z$_2$-graded quantum plane or the quantum superplane 
is defined as an associative algebra whose even coordinate $x$ and the odd 
(Grassmann) coordinate $\theta$ satisfy 
$$ x \theta = q \theta x, \qquad \theta^2 = 0 $$
where $q$ is a nonzero complex deformation parameter. 

One of the possible ways to generalize the quantum superplane is to increase 
the power of nilpotency of its odd generator. So, a possible generalization 
can be defined as an associative unital algebra generated by $x$ and $\theta$ 
satisfying 
$$ x \theta = q \theta x, \qquad \theta^3 = 0. \eqno(1)$$
Here, the coordinate $x$ with respect to the Z$_3$-grading is of grade 0 and 
the coordinate $\theta$ with respect to the Z$_3$-grading is of grade 1. 

The quantum superplane underlies a noncommutative differential calculus on a 
smooth manifold with exterior differential {\sf d} satisfying ${\sf d}^2 = 0.$ 
So the above mentioned generalization of the superplane raise the natural 
question of possible generalization of differential calculus to one with 
exterior differential {\sf d} satisfying ${\sf d}^3 = 0$. From an algebraic 
point of view, a sufficent algebraic structure underlying a differential 
calculus is the notion of the Z$_3$-graded differential algebra. Therefore 
we can generalize the differential calculus with the help of an appropriate 
generalization of Z$_3$-graded differential algebra. 

Elementary properties of Z$_2$-graded quantum superplane are described in Ref. 
13. We state briefly the properties we are going to need in this work. 

Let ${\cal A}$ be a free unital associative algebra generated by two elements 
$x$, $\theta$ obeying the relations (1). We know that the algebra ${\cal A}$ 
is a graded Hopf algebra with the following co-structures [13]: the coproduct 
$\Delta: {\cal A} \longrightarrow {\cal A} \otimes {\cal A}$ 
is defined by 
$$\Delta(x) = x \otimes x, \qquad 
  \Delta(\theta) = \theta \otimes x + x \otimes \theta, \eqno(2)$$
$$\Delta(1) = 1 \otimes 1. $$
The counit $\epsilon: {\cal A} \longrightarrow {\cal C}$ is given by 
$$\epsilon(x) = 1, \qquad \epsilon(\theta) = 0. \eqno(3)$$
We extend the algebra ${\cal A}$ by including inverse of $x$ which obeys 
$$x x^{-1} = 1 = x^{-1} x. $$
If we extend the algebra ${\cal A}$ by adding the inverse of $x$ 
then the algebra ${\cal A}$ admits a coinverse 
$S: {\cal A} \longrightarrow {\cal A}$ defined by 
$$S(x) = x^{-1}, \qquad S(\theta) = - x^{-1} \theta x^{-1}. \eqno(4)$$
Note that 
$$\Delta(x^{-1}) = x^{-1} \otimes x^{-1}.$$
It is not difficult to verify the following properties of co-structures: 
$$(\Delta \otimes \mbox{id}) \circ \Delta = 
  (\mbox{id} \otimes \Delta) \circ \Delta, $$
$$\mu \circ (\epsilon \otimes \mbox{id}) \circ \Delta 
  = \mu' \circ (\mbox{id} \otimes \epsilon) \circ \Delta, \eqno(5)$$
$$m \circ (S \otimes \mbox{id}) \circ \Delta = \epsilon 
  = m \circ (\mbox{id} \otimes S) \circ \Delta $$
where id denotes the identity mapping, 
$$\mu : {\cal C} \otimes {\cal A} \longrightarrow {\cal A}, \qquad 
  \mu' : {\cal A} \otimes {\cal C} \longrightarrow {\cal A} $$
are the canonical isomorphisms, defined by 
$$\mu(c \otimes a) = c a = \mu'(a \otimes c), \qquad \forall a \in {\cal A}, 
  \quad \forall c \in {\cal C} $$
and $m$ is the multiplication map 
$$m : {\cal A} \otimes {\cal A} \longrightarrow {\cal A}, \qquad 
  m(a \otimes b) = ab. \eqno(6)$$
The multiplication in ${\cal A} \otimes {\cal A}$ is defined with the rule 
$$(A \otimes B) (C \otimes D) = j^{grad(B) grad(C)} AC \otimes BD. 
  \eqno(7)$$

\vspace*{0.3cm}\noindent
{\bf 3. Differential calculi on the Z$_3$-graded quantum superplane }

\noindent
In this section, we shall build up the noncommutative differential calculus 
on the Z$_3$-graded quantum superplane. This involves functions on the 
superplane, differentials and differential forms. So we have to define a 
linear operator {\sf d} which acts on the functions of the coordinates of the 
Z$_3$-graded quantum superplane. For the definition, it is sufficent to define 
the action of {\sf d} on the coordinates and on their products. 

We postulate that the linear operator {\sf d} applied to $x$ produces a 1-form 
whose Z$_3$-grade is 1, by definition. Similarly, application of {\sf d} to 
$\theta$ produces a 1-form whose Z$_3$-grade is 2. We shall denote the 
obtained quantities by ${\sf d} x$ and ${\sf d} \theta$, respectively. When 
the linear operator {\sf d} is applied to ${\sf d} x$ (or twice by iteration 
to $x$) it will produce a new entity which we shall call a 1-form of grade 2, 
denoted by ${\sf d}^2 x$ and to ${\sf d} \theta$ produces a 1-form of grade 0, 
modulo 3, denoted by ${\sf d}^2 \theta$. Finally, we require that 
${\sf d}^3 = 0$. 

\vspace*{0.3cm}\noindent
{\bf 3. 1 Differential algebra} 

\noindent
Let us begin the ordering the properties of the exterior differential. The 
exterior differential {\sf d} is an operator which gives the mapping from the 
generators of the Z$_3$-graded quantum superplane to the differentials 
$${\sf d} : a \mapsto {\sf d} a, \qquad a \in \{x,\theta\}. $$
We demand that the exterior differential {\sf d} has to satisfy two 
properties:  
$${\sf d}^3 = 0 \eqno(8)$$
and the Z$_3$-graded Leibniz rule 
$${\sf d}(f g) = ({\sf d} f) g + j^{grad(f)} f ({\sf d} g). \eqno(9)$$

It is well known that in classical differential calculus, functions commute 
with differentials. From an algebraic point of view, the space of 1-forms is 
a free finite bimodule over the algebra of smooth functions generated by the 
first order differentials and the commutativity shows how its left and right 
structure are related to each other. 

In order to establish a noncommutative differential calculus on the 
Z$_3$-graded quantum superplane, we assume that the commutation relations 
between the coordinates and their differentials are in the following form: 
$$x ~{\sf d}x = X {\sf d}x ~x, $$
$$x ~{\sf d}\theta = A {\sf d}\theta ~x + B {\sf d}x ~\theta, $$
$$\theta ~{\sf d}x = C {\sf d}x ~\theta + D {\sf d}\theta ~x, $$
$$\theta ~{\sf d}\theta = Y {\sf d}\theta ~\theta. \eqno(10)$$

The coefficients $A$, $B$ and $F_{ik}$ will be determined in terms of the 
complex deformation parameter $q$ and $j$. To find them we shall use the 
covariance of the noncommutative differential calculus. 

Since we assume that ${\sf d}^3 = 0$ and ${\sf d}^2 \neq 0$, in order to 
construct a self-consistent theory of differential forms it is necessary to 
add to the first order differentials of coordinates ${\sf d} x$, 
${\sf d} \theta$ a set of {\it second order differentials} ${\sf d}^2 x$, 
${\sf d}^2 \theta$. Apperance of higher order differentials is a peculiar 
property of a proposed generalization of differential forms. This has as a 
consequence certain problems. 

Now, we assume that {\sf d} is no more the classical exterior differential, 
i.e. ${\sf d}^2 \neq 0$. For example, if we take a particular 1-form 
$\theta {\sf d} x$ and apply to it the exterior differential {\sf d}, we 
obtain 
$${\sf d} (\theta ~{\sf d} x) = {\sf d} \theta ~{\sf d} x + 
  j \theta ~{\sf d}^2 x.$$
Therefore, differentiating (10) with regard to the Z$_3$-graded Leibniz rule 
(9) one gets 
$$x ~{\sf d}^2 x = X ~{\sf d}^2 x ~x + (jX - 1) ({\sf d} x)^2,$$
$$x ~{\sf d}^2 \theta = A {\sf d}^2 \theta ~x + B {\sf d}^2 x ~\theta 
  + (j^2 A + j B F - F) {\sf d} \theta ~{\sf d} x, $$
$$\theta ~{\sf d}^2 x = j^{-1} C {\sf d}^2 x ~\theta + 
  j^{-1} D {\sf d}^2 \theta ~x + 
  (jD + CF - j^{-1}) {\sf d}\theta ~{\sf d} x, $$ 
$$\theta ~{\sf d}^2 \theta = Y j^{-1} {\sf d}^2 \theta ~\theta + 
  (jY - j^{-1}) ({\sf d} \theta)^2. \eqno(11\mbox{a}) $$
Here, we have assumed that 
$${\sf d} x ~{\sf d} \theta = F ~{\sf d} \theta ~{\sf d} x, 
  \qquad ({\sf d} x)^3 = 0, \eqno(11\mbox{b})$$
where $F$ is a parameter that shall described later. 

The relations (11a) are {\it not homogeneous} in the sense that the 
commutation relations between the coordinates and second order differentials 
include first order differentials as well. Later, we shall see that the 
commutation relations between the coordinates and their second order 
differentials can be made homogeneous. They will not include first order 
differentials by removing them using the covariance of the noncommutative 
differential calculus. 

We now return to the relations (11a). Applying the exterior differential 
{\sf d} to the relations (11a), we obtain 
$${\sf d} x ~{\sf d}^2 x = j^{-2} {\sf d}^2 x ~{\sf d} x, $$
$${\sf d}x ~{\sf d}^2 \theta = - {A\over Q} ~{\sf d}^2 \theta ~{\sf d}x + 
  {{1 + B - j^2 A F^{-1}}\over Q} ~{\sf d}^2 x ~{\sf d} \theta, $$
$${\sf d} \theta ~{\sf d}^2 x = - {C\over {Q'}} ~{\sf d}^2 x ~{\sf d} \theta 
   + {{D - CF + j^{-1}}\over {Q'}} ~{\sf d}^2 \theta ~{\sf d} x, $$ 
$${\sf d} \theta ~{\sf d}^2 \theta = {\sf d}^2 \theta ~{\sf d} \theta, 
  \eqno(11\mbox{c})$$
where 
$$Q = AF^{-1} + j^2 (1 + B) \quad \mbox{and} \quad Q' = D + j^2 (1 + CF).$$
As a consequence, the second order differentials have to satisfy the following 
relation 
$${\sf d}^2 x ~{\sf d}^2 \theta = j F {\sf d}^2 \theta {\sf d}^2 x. 
  \eqno(11\mbox{d})$$

\vspace*{0.3cm}\noindent
{\bf 3. 2 Covariance}

\noindent
We see from the above relations (11a) that the commutation relations between 
the generators of ${\cal A}$ and their second order differentials are 
{\it not homogeneous} in the sense that they include first order 
differentials. In order to homogenize the relations (11a), we shall consider 
the covariance of the noncommutative differential calculus. 

We first note that consistency of a differential calculus with commutation 
relations (1) means that the differential algebra is a graded associative 
algebra generated by the elements of the set 
$\{x, \theta, {\sf d} x, {\sf d} \theta, {\sf d}^2 x, {\sf d}^2 \theta\}$. 

Let $\Omega({\cal A})$ be a free left module over the algebra ${\cal A}$ 
generated by the elements of the set 
$\{x, \theta, {\sf d} x, {\sf d} \theta, {\sf d}^2 x, {\sf d}^2 \theta\}$. The 
module $\Omega({\cal A})$ becomes a unital associative algebra if one defines 
a multiplication law on $\Omega({\cal A})$ by the relations (1) and (11). 

We consider a map 
$\phi_L : \Omega({\cal A}) \longrightarrow {\cal A} \otimes \Omega({\cal A})$ 
such that 
$$\phi_L \circ {\sf d} = (\tau \otimes {\sf d}) \circ \Delta, \eqno(12)$$
where $\tau: \Omega({\cal A}) \longrightarrow \Omega({\cal A})$ is the linear 
map of degree zero which gives 
$$\tau(a) = j^{grad(a)} a, \qquad \forall a \in \Omega({\cal A}). \eqno(13)$$
Thus we have 
$$\phi_L({\sf d} x) = x \otimes {\sf d} x, \qquad 
  \phi_L({\sf d}\theta) = j \theta \otimes {\sf d} x + x \otimes {\sf d}\theta.
  \eqno(14)$$
We now define a map $\Delta_L$ as follows: 
$$\Delta_L(a_1 ~{\sf d}b_1 + {\sf d}b_2 ~a_2) = 
  \Delta(a_1) \phi_L({\sf d}b_1) + \phi_L({\sf d}b_2) \Delta(a_2). \eqno(15)$$

We now apply the linear map $\Delta_L$ to relations (15): 
$$\Delta_L(x ~{\sf d} x) = \Delta(x) \phi_L ({\sf d} x) = 
  X \phi_L ({\sf d} x) \Delta(x) = X \Delta_L ({\sf d} x ~x), $$
\begin{eqnarray*}
\Delta_L(x ~{\sf d} \theta) 
& = & \Delta(x) \phi_L ({\sf d} \theta) \\
& = & A \phi_L ({\sf d} \theta) \Delta(x) + 
      B \phi_L({\sf d} x) \Delta(\theta) + 
      j (X - q^{-1} A - B) x \theta \otimes {\sf d}x ~x \\
& = & A \Delta_L ({\sf d} \theta ~x) + B \Delta_L({\sf d} x ~\theta) + 
      j (X - q^{-1} A - B) x \theta \otimes {\sf d}x ~x,
\end{eqnarray*}
\begin{eqnarray*}
\Delta_L(\theta ~{\sf d} x) 
& = & \Delta(\theta) \phi_L ({\sf d} x) \\
& = & C \phi_L ({\sf d} x) \Delta(\theta) + D \phi_L ({\sf d}\theta) \Delta(x) 
      + (X - qjC - jD) \theta x \otimes {\sf d}x ~x \\
& = & C \Delta_L({\sf d} x ~\theta) + D \Delta_L ({\sf d} \theta ~x) 
      + (X - qjC - jD) \theta x \otimes {\sf d}x ~x 
\end{eqnarray*}
and 
\begin{eqnarray*}
\Delta_L(\theta {\sf d} \theta) 
& = & \Delta(\theta) \phi_L ({\sf d} \theta) \\
& = & Y \phi_L ({\sf d} \theta) \Delta(\theta) + 
      (B - j Y + q j^2 C) \theta x \otimes {\sf d} x ~\theta\\
&   & + j (X - jY) \theta^2 \otimes {\sf d} x ~x + 
      (A + qj^2 D - qj^2 Y) x \theta \otimes {\sf d}\theta ~x \\
& = & Y \Delta_L({\sf d} \theta ~\theta) + 
      (B - j Y + q j^2 C) \theta x \otimes {\sf d} x ~\theta\\
&   & + j (X - jY) \theta^2 \otimes {\sf d} x ~x + 
      (A + qj^2 D - qj^2 Y) x \theta \otimes {\sf d}\theta ~x. 
\end{eqnarray*}
We see from the last three relations that in order to have left covariance 
$D$ must be zero. Then with $Y$ arbitrary 
$$A = j^2 q Y, \qquad C = q^{-1} Y, $$
$$B  = j (1 - j) Y, \qquad X = j Y. \eqno(16)$$
On the other hand, the action of {\sf d} on $\theta^3 = 0$ gives 
$$1 + j Y + j ^2 Y^2 = 0. \eqno(17a)$$
So, 
$$Y = j \qquad \mbox{or} \qquad Y = j^2. \eqno(17b)$$
For $Y = j^2$, the relations (11a) are not homogeneous. Hence, we must take 
$Y = j$. 

Also, since 
\begin{eqnarray*}
\Delta_L({\sf d} x {\sf d} \theta) 
& = & \Delta_L({\sf d} x) \Delta_L ({\sf d} \theta) \\
& = & F \Delta_L ({\sf d} \theta) \Delta_L({\sf d} x) + 
      j (q j - F) \theta x \otimes ({\sf d}x)^2, 
\end{eqnarray*}
we must have 
$$F - q j = 0. \eqno(18)$$
Here, we used that 
$$(x \otimes {\sf d} x) (\theta \otimes {\sf d} x) = 
  j x \theta \otimes ({\sf d} x)^2.$$

The relations (10) and (11) are explicity as follows: the commutation 
relations of variables and their differentials are 
$$x~ {\sf d}x = j^2 ~{\sf d}x~ x, $$
$$x~ {\sf d}\theta = q ~{\sf d}\theta~ x + (j^2 - 1) {\sf d}x~ \theta, $$
$$\theta~ {\sf d}x = j q^{-1} {\sf d}x~ \theta, $$ 
$$\theta~ {\sf d}\theta = j ~{\sf d}\theta~ \theta, \eqno(19)$$
and among those first order differentials are 
$${\sf d} x ~{\sf d} \theta = j q ~{\sf d} \theta ~{\sf d} x, 
  \qquad ({\sf d} x)^3 = 0, \eqno(20)$$

The commutation relations between variables and second order differentials are 
$$x ~{\sf d}^2 x = j^2 {\sf d}^2 x ~x,$$
$$x ~{\sf d}^2 \theta = q{\sf d}^2 \theta ~x + (j^2 - 1){\sf d}^2 x ~\theta, $$
$$\theta ~{\sf d}^2 x = q^{-1} {\sf d}^2 x ~\theta, $$ 
$$\theta ~{\sf d}^2 \theta = {\sf d}^2 \theta ~\theta. \eqno(21) $$

The commutation relations between first order and second order differentials are 
$${\sf d} x ~{\sf d}^2 x = j^{-2} {\sf d}^2 x ~{\sf d} x, $$
$${\sf d}x ~{\sf d}^2 \theta = q ~{\sf d}^2 \theta ~{\sf d}x + 
  (j - j^{-1}) ~{\sf d}^2 x ~{\sf d} \theta, $$
$${\sf d} \theta ~{\sf d}^2 x = j^2 q^{-1} {\sf d}^2 x ~{\sf d} \theta, $$ 
$${\sf d} \theta ~{\sf d}^2 \theta = {\sf d}^2 \theta ~{\sf d} \theta, 
  \eqno(22)$$
and those among the second order differentials are 
$${\sf d}^2 x ~{\sf d}^2 \theta = j^2 q~ {\sf d}^2 \theta ~{\sf d}^2 x. 
  \eqno(23)$$

Now, it can be checked that the linear map $\Delta_L$ leaves invariant the 
relations (19)-(23). One can also check that the following identities 
are satisfied: 
$$(\mbox{id} \otimes \Delta_L) \circ \Delta_L = 
  (\Delta \otimes \mbox{id}) \circ \Delta_L, \qquad 
m \circ (\epsilon \otimes \mbox{id}) \circ \Delta_L = \mbox{id}. \eqno(24)$$
We call as left coaction the map $\Delta_L$. The map $\Delta_L$ makes the 
$\Omega({\cal A})$ a left ${\cal A}$-module. 
So, the pair $(\Omega({\cal A}), \Delta_L)$ is a left-covariant left 
${\cal A}$-module over Hopf algebra ${\cal A}$. Hovewer the pair 
$(\Omega({\cal A}), {\sf d})$ is a differential calculus over ${\cal A}$, and 
{\sf d} is a left comodule map, i.e. for all $a \in {\cal A}$, 
$$(\tau \otimes {\sf d})\circ \Delta(a) = \Delta_L({\sf d} a). \eqno(25)$$
Consequently, the triple $(\Omega({\cal A}), {\sf d}, \Delta_L)$ is a left 
covariant differential calculus over the Hopf algebra ${\cal A}$. 

\vspace*{0.3cm}\noindent
{\bf 4. Cartan-Maurer one-forms on ${\cal A}$}

\noindent
In this section we shall define two forms using the generators of 
$\cal A$ and show that they are left-invariant. 
If we call them $w$ and $u$ then one can define them as follows [13]: 
$$w = {\sf d}x ~x^{-1}, \qquad 
  u = {\sf d}\theta ~x^{-1} - {\sf d}x ~x^{-1} \theta x^{-1}. \eqno(26)$$
The elements $w$ and $u$ with the generators of $\cal A$ satisfy the following 
rules 
$$x w = j^2 w x, \qquad \theta w = j w \theta, $$
$$x u = q u x, \qquad \theta u = j q u \theta. \eqno(27)$$
The first order differentials with 1-forms satisfy the following relations 
$$w ~{\sf d} x = j {\sf d} x ~w, \qquad u ~{\sf d} x = q^{-1} {\sf d} x ~u,$$
$$w ~{\sf d} \theta = j^2 {\sf d} \theta ~w + q^{-1} (1 - j) {\sf d} x ~u,$$
$$u ~{\sf d} \theta = q^{-1} {\sf d} \theta ~u + 
  q^{-2} (1 - j) {\sf d} x ~u ~\theta x^{-1} \eqno(28\mbox{a}),$$
and with second order differentials 
$$w ~{\sf d}^2 x = j^2 {\sf d}^2 x ~w, \qquad 
  u ~{\sf d}^2 x = q^{-1} {\sf d}^2 x ~u,$$
$$w ~{\sf d}^2 \theta = {\sf d}^2 \theta ~w + 
  q^{-1} (j - j^{-1}) {\sf d}^2 x ~u,$$
$$u ~{\sf d}^2 \theta = q^{-1} {\sf d}^2 \theta ~u + 
  q^{-2} (1 - j) {\sf d}^2 x ~u ~\theta x^{-1}. \eqno(28\mbox{b})$$
The commutation rules of the elements $w$ and $u$ are 
$$w^3 = 0, \qquad wu = uw. \eqno(29)$$

The elements $w$ and $u$ are both left invariant with the following 
structures: 
$$\Delta_L(w) = 1 \otimes w, \qquad \Delta_L(u) = 1 \otimes u. \eqno(30) $$
The counit $\epsilon$ is given by [13]
$$\epsilon(w) = 0, \qquad \epsilon(u) = 0 \eqno(31)$$
and the coinverse $S$ is defined by 
$$S(w) = - w, \qquad S(u) = - u. \eqno(32)$$
One can easily to check that the following properties are satisfied: 
$$(\mbox{id} \otimes \Delta_L) \circ \Delta_L = (\Delta \otimes \mbox{id}) 
  \circ \Delta_L, $$
$$ m \circ (\epsilon \otimes \mbox{id}) \circ \Delta_L = \mbox{id},$$
$$m \circ (S \otimes \mbox{id}) \circ \Delta_L = \mbox{id}. $$
Note that the commutation relations (27)-(29) are compatible with 
$\Delta_L$, $\epsilon$ and $S$, in the sense that 
$\Delta_L(x w) = \Delta(x) \Delta_L(w) = j^2 \Delta_L(w x)$, 
$\Delta_L(w^3) = 0$ and so on. 

%\vfill\eject
\vspace*{0.3cm}\noindent
{\bf 5. Quantum Lie superalgebra} 

\noindent
The commutation relations of Cartan-Maurer forms allow us to construct the 
algebra of the generators. In order to obtain the quantum Lie superalgebra 
of the algebra generators we first write the Cartan-Maurer forms as 
$${\sf d} x = w x, \qquad {\sf d} \theta = w \theta + u x. \eqno(33)$$
The differential {\sf d} can then the expressed in the form 
$${\sf d} = w T + u \nabla. \eqno(34)$$
Here $T$ and $\nabla$ are the quantum Lie superalgebra generators. 
We now shall obtain the commutation relations of these generators. 
Considering an arbitrary function $f$ of the coordinates of the Z$_3$-graded 
quantum superplane and using that ${\sf d}^3 = 0$ one has 
$${\sf d}^2 f = ({\sf d} w) ~T f + ({\sf d} u) ~\nabla f + 
   j w ~{\sf d} T f + j^2 u ~{\sf d} \nabla f, $$
and 
$${\sf d}^3 f = j^{-1} w ~{\sf d}^2 ~T f + j u ~{\sf d}^2 \nabla f 
   - {\sf d} w ~{\sf d} T f - j {\sf d} u ~{\sf d} \nabla f + 
   {\sf d}^2 w ~T f + {\sf d}^2 u ~\nabla f.$$
So we need the two-forms. Applying the exterior differential {\sf d} to the 
relations (26) one has 
$${\sf d} w = {\sf d}^2 x ~x^{-1} - j w^2, $$
$${\sf d} u = {\sf d}^2 \theta ~x^{-1} - {\sf d}^2 x x^{-1} \theta x^{-1} 
  + u w. \eqno(35)$$
Also, since 
$$w ~{\sf d} w = j {\sf d} w ~w, $$
$$w ~{\sf d} u = j^2 {\sf d} u ~w + (j - j^{-1}) {\sf d} w ~u,$$
$$u ~{\sf d} w = {\sf d} w ~u, \qquad u ~{\sf d} u = {\sf d} u ~u,$$
we have 
$${\sf d}^2 w = 0, \qquad 
 {\sf d}^2 u = 0. \eqno(36)$$
Using the Cartan-Maurer equations we find the following commutation 
relations for the quantum Lie superalgebra: 
$$T \nabla = \nabla T, \qquad \nabla^3 = 0. \eqno(37)$$

The commutation relations (37) of the algebra generators should be consistent 
with monomials of the coordinates of the Z$_3$-graded quantum superplane. 
To do this, we evaluate the commutation relations between the generators of 
algebra and the coordinates. The commuation relations of the generators with 
the coordinates can be extracted from the Z$_3$-graded Leibniz rule: 
\begin{eqnarray*}
{\sf d} (x f) 
& = & ({\sf d} x) f + x ({\sf d} f) \\
& = & w (x + j^2 x T) f + u (q x \nabla) f \\ 
& = & (w T + u \nabla) x f \hspace*{7.8cm}{(38)}
\end{eqnarray*}
and 
\begin{eqnarray*}
{\sf d} (\theta f) 
& = & ({\sf d} \theta) f + j \theta ({\sf d} f) \\
& = & w (\theta + j^2 \theta T) f + u (x + q j^2 \theta \nabla) f \\ 
& = & (w T + u \nabla) \theta f. \hspace*{7.7cm}{(39)}
\end{eqnarray*}
This yields 
$$T x  = x + j^2 x T, \qquad T \theta  = \theta + j^2 \theta T, $$
$$\nabla x = q x \nabla, \qquad \nabla \theta = x + qj^2 \theta \nabla. 
  \eqno(40)$$

We know that the differential operator {\sf d} satisfies the Z$_3$-graded 
Leibniz rule. Therefore, the generators $T$ and $\nabla$ are endowed with a 
natural coproduct. To find them, we need to the following commutation 
relations 
$$T x^m = {{1 - j^{2m}}\over {1 - j^2}} x^m + j^{2m} x^m T, \eqno(41)$$
$$\nabla x^m = q^m x^m \nabla, \eqno(42)$$
where use was made of (40). The relation (41) is understood as an operator 
equation. This implies that when $T$ acts on arbitrary monomials $x^m \theta$, 
$$T (x^m \theta) = {{1 - j^{2m + 2}}\over {1 - j^2}} (x^m \theta) + 
  j^{2m + 2} (x^m \theta) T, \eqno(43)$$
from which we obtain 
$$T = {{1 - j^{2N}}\over {1 - j^2}},  \eqno(44)$$
where $N$ is a number operator acting on a monomial as 
$$N(x^m \theta) = (m + 1) x^m \theta. \eqno(45)$$
We also have 
$$\nabla(x^m \theta) = q^m x^{m+1} + j^2 q^{m+1}(x^m \theta)\nabla. \eqno(46)$$

So, applying the Z$_3$-graded Leibniz rule to the product of functions $f$ and 
$g$, we write 
$${\sf d} (f g) = [(w T + u \nabla) f] g + j^{grad(f)} f (w T + u \nabla) g 
  \eqno(47)$$
with help of (34). From the commutation relations of the Cartan-Maurer forms 
with the coordinates of the Z$_3$-graded quantum superplane, we can compute 
the corresponding relations of $w$ and $u$ with functions of the coordinates. 
From (27) we have 
$$f w = j^{2N - 1} w f, \qquad f u = j q^N u f \eqno(48)$$
where $f = x^m \theta$. Inserting (48) 
in (47) and equating coefficients of the Cartan-Maurer forms, we get 
$$T(f g) = (T f) g + j^{grad(f)} j^{2N - 1} f (T g), $$
$$\nabla(f g) = (\nabla f) g + j^{grad(f)} j q^N f (\nabla g). \eqno(49)$$
Consequently, we have the coproduct 
$$\Delta(T) = T \otimes 1 + j^{- N} \otimes T, $$
$$\Delta(\nabla) = \nabla \otimes 1 + j^2 q^N \otimes \nabla. \eqno(50)$$

\vspace*{0.3cm}\noindent
{\bf 6. Conclusion}

\noindent
To conclude, we introduce here commutation relations between the coordinates 
of the Z$_3$-graded quantum superplane and their partial derivatives and thus 
illustrate the connection between the relations in section 5, and the 
relations which will be now obtained. 

To proceed, let us obtain the relations of the coordinates with their partial 
derivatives. We know that the exterior differential {\sf d} can be expressed 
in the form 
$${\sf d} f = ({\sf d} x ~\partial_x + {\sf d} \theta ~\partial_\theta) f. 
  \eqno(51)$$
Then, for example, 
\begin{eqnarray*}
{\sf d} (x f) 
& = & {\sf d} x ~f + x ~{\sf d} f\\
& = & {\sf d} x ~[1 + j^2 x \partial_x + (j^2 - 1) \theta \partial_\theta] f 
      + q {\sf d} \theta ~x \partial_\theta f \\
& = & ({\sf d} x ~\partial_x x + {\sf d} \theta ~\partial_\theta x) f
\end{eqnarray*}
so that 
$$\partial_x x = 1 + j^2 x \partial_x + (j^2 - 1) \theta \partial_\theta,
  \qquad \partial_x \theta = j^2 q^{-1} \theta \partial_x,$$
$$\partial_\theta x = q x \partial_\theta, \qquad 
  \partial_\theta \theta = 1 + j^2 \theta \partial_\theta. \eqno(52)$$
The commutation relations between derivatives are 
$$\partial_x \partial_\theta = j q \partial_\theta \partial_x, \qquad 
  \partial_\theta^3 = 0. \eqno(53)$$
The Z$_3$-graded Hopf algebra structure for $\partial$ is given by 
$$\Delta(\partial_x) = \partial_x \otimes \partial_x, \qquad 
  \Delta(\partial_\theta) = \partial_\theta \otimes \partial_x + 
  \partial_x \otimes \partial_\theta, $$
$$\epsilon(\partial_x) = 1, \qquad \epsilon(\partial_\theta) = 0, \eqno(54)$$ 
$$S(\partial_x) = \partial_x^{-1}, \qquad 
  S(\partial_\theta) = - \partial_x^{-1} \partial_\theta \partial_x^{-1}, $$
provided that the formal inverse $\partial_x^{-1}$ exists. 
However these co-maps do not leave invariant the relations (52). 

We know, from section 5, that the exterior differential {\sf d} can be 
expressed in the form (34), which we repeat here, 
$${\sf d} f = (w T + u \nabla) f. \eqno(55)$$
Considering (51) together (55) and using (33) one has 
$$T = x \partial_x + \theta \partial_\theta, \qquad 
  \nabla = x \partial_\theta. \eqno(56)$$
Using the relations (52) and (53) one can check that the relations of the 
generators in (56) coincide with (37). It can also be verified that, the 
action of the generators in (56) on the coordinates coincide with (40). 

The Z$_3$-graded noncommutative differential geometry we have constructed 
satisfies all expectations for such a structure. In particular all Hopf 
algebra axioms are satisfied without any modification. Moreover the extension 
of the structure presented in this paper can be generalized to Z$_N$. 

\vspace*{0.3cm}\noindent
{\bf Appendix: Quantum matrices in Z$_3$-graded space} 

\noindent
In this appendix we shall investigate the quantum supermatrices in 
Z$_3$-graded quantum superplane. We know, from section 2, that the 
Z$_3$-graded quantum superplane is generated by coordinates $x$ and $\theta$, 
and the commutation rules (1), which we repeat here, 
$$x \theta = q \theta x, \qquad \theta^3 = 0. \eqno(57)$$
These relations define a deformation of the algebra of functions on the 
superplane generated by $x$ and $\theta$, and we have denoted it by $\cal A$. 
The dual Z$_3$-graded quantum superplane $\cal A^\star$ is generated by 
$\varphi$ and $y$ with the relations 
$$\varphi y = q j y \varphi, \qquad \varphi^3 = 0, \eqno(58)$$
where ${\sf d} x = \varphi$ and ${\sf d} \theta = y$ in (20). 

Let $T$ be a 2x2 (super)matrix in Z$_3$-graded space, 
$$T = \left(\matrix{a & \beta \cr \gamma & d}\right) \eqno(59)$$
where $a$ and $d$ with respect to the Z$_3$-grading are of grade 0, and 
$\beta$ and $\gamma$ with respect to the Z$_3$-grading are of grade 2 and 
of grade 1, respectively. We now consider linear transformations with the 
following properties: 
$$T : {\cal A} \longrightarrow {\cal A}, \qquad 
  T : {\cal A^\star} \longrightarrow {\cal A^\star} \eqno(60)$$
The action on the elements of $\cal A$ of $T$ is 
$\left(\matrix{x' \cr \theta'}\right) = 
  \left(\matrix{a & \beta \cr \gamma & d}\right) 
  \left(\matrix{x \cr \theta}\right).$ 
We assume that the entries of $T$ are $j$-commutative with the elements of 
$\cal A$, i.e. for example, 
$$a x = x a, \qquad \theta \beta = j^2 \beta \theta, $$
etc. As a consequence of the linear transformations in (60) the elements 
$$x' = a x + \beta \theta, \qquad \theta' = \gamma x + d \theta \eqno(61)$$
should satisfy the relations (57). Applying the exterior differential {\sf d} 
to the relation (61) one has 
$$\varphi' = a\varphi + j^2\beta y, \qquad y' = j\gamma\varphi+d y. \eqno(62)$$
These elements must satisfy the relations (58). Consequently, we have 
the following commutation relations between the matrix elements of $T$: 
$$a \beta = j^{-1} q^{-1} \beta a, \qquad d \beta = j q^{-1} \beta d, $$
$$a \gamma = q \gamma a, \qquad d \gamma = q \gamma d,$$
$$a d = d a + q^{-1} (1 - j) \beta \gamma, \qquad 
  \beta \gamma = q^2 \gamma \beta, \qquad \gamma^3 = 0. \eqno(63)$$
We shall denote with GL$_{q,j}(1\vert 1)$ the quantum supergroup in 
Z$_3$-graded space determined by generators $a$, $\beta$, $\gamma$, $d$ 
satisfying the commutation relations (63). 

Note that these relations can be obtained from the requirement that $\cal A$ 
and ${\cal A^\star}$ have to be covariant under the left coactions 
$$\delta : {\cal A} \longrightarrow \mbox{GL}_{q,j}(1\vert 1) \otimes {\cal A},
  \qquad \delta^\star : {\cal A^\star} \longrightarrow 
  \mbox{GL}_{q,j}(1\vert 1) \otimes {\cal A^\star} \eqno(64)$$
such that 
$$\delta(x) = a \otimes x + \beta \otimes \theta, \qquad 
  \delta(\varphi) = a \otimes \varphi + j^2 \beta \otimes y, $$
$$\delta(\theta) = \gamma \otimes x + d \otimes \theta, \qquad 
  \delta(y) = j \gamma \otimes \varphi + d \otimes y, \eqno(65)$$
provided that the entries of $T$ are $j$-commutative with the elements of 
${\cal A}$ and ${\cal A^\star}$. 

Note that the relations (63) are slightly different from the results of 
Ref. 12. The reason for this difference is that in Ref. 12, 
since he assumed that commutation relations of 
the differentials are 
$$ {\sf d} x {\sf d} \theta = r^{-1} {\sf d} \theta {\sf d} x, \qquad 
  ({\sf d} x)^2 = 0 =({\sf d} \theta)^2, $$
the commutation relations among the matrix elements of a matrix in 
Z$_3$-graded space were obtained via the use of them. On the other hand, we 
use the commutation relations of the coordinates of the Z$_3$-graded quantum 
superplane with their differentials. 

An interesting problem is the construction of a differential calculus on the 
Z$_3$-graded quantum supergroup GL$_{q,j}(1\vert 1)$ using the methods of 
this paper and Ref. 15. Work on this issue is in progress. 

\noindent
{\bf Acknowledgement}

\noindent
This work was supported in part by {\bf T. B. T. A. K.} the 
Turkish Scientific and Technical Research Council.

\end{document}